\documentclass{article}
\usepackage{latexsym}
\usepackage{epsfig}
\usepackage{amsmath}
\usepackage{amssymb}
\usepackage{amsthm}
\usepackage{eufrak}
\newtheorem{theorem}{Theorem}[section]

\newtheorem{lemma}[theorem]{Lemma}

\theoremstyle{remark}

\theoremstyle{definition}

\DeclareMathOperator\Supp{supp}

\begin{document}

\title{Extremal copositive matrices with minimal zero supports of cardinality two}

\author{Roland Hildebrand \thanks{%
LJK / CNRS, B\^atiment IMAG, 700 avenue Centrale, Domaine Universitaire, 38041 Saint-Martin-d'H\`eres, France
({\tt roland.hildebrand@univ-grenoble-alpes.fr}).}}

\maketitle

\begin{abstract}
Let $A \in {\cal C}^n$ be an extremal copositive matrix with unit diagonal. Then the minimal zeros of $A$ all have supports of cardinality two if and only if the elements of $A$ are all from the set $\{-1,0,1\}$. Thus the extremal copositive matrices with minimal zero supports of cardinality two are exactly those matrices which can be obtained by diagonal scaling from the extremal $\{-1,0,1\}$ unit diagonal matrices characterized by Hoffman and Pereira in 1973.
\end{abstract}



\section{Introduction}

An element $A$ of the space ${\cal S}^n$ of real symmetric $n \times n$ matrices is called \emph{copositive} if $x^TAx \geq 0$ for all vectors $x \in \mathbb R_+^n$. The set of such matrices forms the \emph{copositive cone} ${\cal C}^n$. This cone plays an important role in non-convex optimization, as many difficult optimization problems can be reformulated as conic programs over ${\cal C}^n$. For a detailed survey of the applications of this cone see, e.g.,~\cite{Duer10,Bomze12}.

Verifying copositivity of a given matrix is a co-NP-complete problem~\cite{MurtyKabadi87}, and the complexity of the copositive cone quickly grows with dimension. In this note we focus on the extreme rays of ${\cal C}^n$. This topic has attracted particular interest already very early in the study of copositive matrices, and a number of families of extreme copositive matrices have been constructed~\cite{HallNewman63,Baumert67,Baston69,HaynsworthHoffman69,HoffmanPereira73}. An element $x$ of a regular convex cone $K$ is called an {\it extremal} element if a decomposition $x = x_1 + x_2$ of $x$ into elements $x_1,x_2 \in K$ is only possible if $x_1 = \lambda x$, $x_2 = (1-\lambda)x$ for some $\lambda \in [0,1]$. The set of positive multiples of an extremal element is called an {\it extreme ray} of $K$. The set of extreme rays is an important characteristic of a convex cone. Its structure, first of all its stratification into a union of manifolds of different dimension, yields much information about the shape of the cone. The extreme rays of a convex cone which is algorithmically difficult to access are especially important if one wishes to check the tightness of inner convex approximations of the cone. Namely, an inner approximation is exact if and only if it contains all extreme rays. Since the extreme rays of a cone determine the facets of its dual cone, they are also important tools for the study of this dual cone. The extreme rays of the copositive cone have been used in a number of papers on its dual, the completely positive cone~\cite{Dickinson11,SMBJS13,BomzeSchachingerUllrich14,BomzeSchachingerUllrich15,SMBBJS15,ShakedMonderer17}.

A important tool in the study of extremal copositive matrices are its zeros~\cite{Diananda62,Baumert66}. A \emph{zero} $u$ of a copositive matrix $A$ is a non-zero nonnegative vector such that $u^TAu = 0$. The \emph{support} $\Supp u$ of a zero $u = (u_1,\dots,u_n)^T \in \mathbb R_+^n$ is the subset of indices $j \in \{1,\dots,n\}$ such that $u_j > 0$. In \cite{Hildebrand14a} we introduced a refined tool, the minimal zeros. Here a zero $u$ of a copositive matrix $A$ is called \emph{minimal} if there exists no zero $v$ of $A$ such that $\Supp v \subset \Supp u$ holds strictly. Up to multiplication by a positive constant there exists only a finite number of minimal zeros for a given copositive matrix. The set of supports of all minimal zeros of a copositive matrix is an informative characteristic of the matrix. In particular, this combinatorial characteristic can assist the classification of the extremal elements of ${\cal C}^n$.

In this note we make a step in this direction by describing the extremal copositive matrices whose minimal zero supports have all cardinality two. We show that every such matrix can be transformed by an automorphism of the copositive cone to a matrix whose elements are from the set $\{-1,0,1\}$. However, the extremal copositive matrices with the latter property have been already completely classified in~\cite{HoffmanPereira73}. This yields a complete description of the former class of extremal copositive matrices. We essentially use the following recent result giving a necessary and sufficient condition of extremality of copositive matrices in terms of their minimal zeros~\cite[Theorem 17]{DickinsonHildebrand16}.

\begin{theorem} \label{th_extremal}
Let $A \in {\cal C}^n$ be a copositive matrix, and let $u^1,\dots,u^m$ be all its minimal zeros, up to multiplication of the zero by a positive constant. Consider the following linear homogeneous system of equations on the matrix $X \in {\cal S}^n$:
\begin{equation} \label{extr_system}
(Xu^j)_k = 0\qquad \forall\ k = 1,\dots,n,\ j = 1,\dots,m\ {\mbox such that}\ (Au^j)_k = 0.
\end{equation}
Then $A$ is an extremal copositive matrix if and only if the solution space of system \eqref{extr_system} has dimension 1. \qed
\end{theorem}

\section{Main result}

We show the two directions of the announced relation between $\{-1,0,1\}$ copositive matrices and matrices with minimal zero supports of cardinality two separately.

{\lemma \label{2to1} Let $A \in {\cal C}^n$ be a copositive matrix with unit diagonal and whose elements are from the set $\{-1,0,1\}$. Then all minimal zeros of $A$ have support of cardinality two. }

\begin{proof}
Since all diagonal elements of $A$ are positive, there exists no zero of $A$ with support of cardinality one.

For the sake of contradiction, suppose $A$ has a minimal zero $u \in \mathbb R_+^n$ with support of cardinality $k > 2$. Without loss of generality, let $\Supp u = \{1,\dots,k\}$. Then all elements of the upper left $k \times k$ submatrix of $A$ equal either 0 or 1. Indeed, suppose there exist $i,j \leq k$ such that $A_{ij} = -1$. Then the sum $v = e^i + e^j$ of the corresponding basis vectors of $\mathbb R^n$ is a zero of $A$ whose support $\{i,j\}$ is a strict subset of $\Supp u$, a contradiction with the minimality of $u$. Hence $A_{ij} \geq 0$ for all $i,j = 1,\dots,k$, and we get
\[ 0 = u^TAu = \sum_{i,j = 1}^k A_{ij}u_iu_j \geq \sum_{i = 1}^k A_{ii}u_i^2 = ||u||_2^2 > 0,
\]
a contradiction.

This completes the proof.
\end{proof}

For the converse direction we shall need the following result~\cite[Corollary 4.4]{DDGH13a}.

{\lemma \label{lem:supp2} Let $A \in {\cal C}^n$ be a copositive matrix with unit diagonal and let $u$ be a zero of $A$ with support $\{i,j\}$. Then $u_i = u_j$. \qed}

\medskip

We are now in a position to prove the following result.

{\lemma \label{1to2} Let $A \in {\cal C}^n$ be an extremal copositive matrix such that all its minimal zeros have support of cardinality two. Then there exists a positive definite diagonal matrix $D$ and an extremal copositive matrix $\Sigma$ with unit diagonal and with all elements in the set $\{-1,0,1\}$ such that $A = D\Sigma D$. }

\begin{proof}
The matrix $A$ has no zeros with support of cardinality one, and therefore all its diagonal elements are positive. Let $D$ be the diagonal matrix with diagonal elements $D_{ii} = \sqrt{A_{ii}}$, $i = 1,\dots,n$, and set $\Sigma = D^{-1}AD^{-1}$.

Note that the linear map given by $X \mapsto D^{-1}XD^{-1}$ is an automorphism of the copositive cone and preserves the property of copositive matrices of being extremal. Hence $\Sigma$ is an extremal copositive matrix with unit diagonal. Note that $u$ is a zero of $A$ if and only if $Du$ is a zero of $\Sigma$, and $\Supp u = \Supp Du$. Hence all minimal zeros of $\Sigma$ have support of cardinality two too. Let $u$ be a minimal zero of $\Sigma$ with support $\{i,j\}$, and let $k \in \{1,\dots,n\}$ be an index. By Lemma \ref{lem:supp2} the equation $(Xu)_k = 0$ on the matrix $X \in {\cal S}^n$ can be written as $X_{ik} + X_{jk} = 0$. Therefore system \eqref{extr_system}, in application to the extremal copositive matrix $\Sigma$, can be written as
\begin{equation} \label{Sigma_sys}
X_{ik} + X_{jk} = 0:\qquad \exists\ \mbox{minimal zero}\ u\ \mbox{of}\ \Sigma\ \mbox{such that}\ \Supp u = \{i,j\},\ (\Sigma u)_k = 0.
\end{equation}

We shall now investigate the solution space of system \eqref{Sigma_sys}. Let $G$ be the graph with the $\frac{n(n+1)}{2}$ independent elements of the real symmetric matrix $X$ as vertices, and with an edge between $X_{ik}$ and $X_{jl}$ if and only if the equation $X_{ik} + X_{jl} = 0$ is among the equations of system \eqref{Sigma_sys}. It is then easily seen that the dimension of the solution space of system \eqref{Sigma_sys} equals the number of bipartite connection components of $G$. Indeed, the elements of $X$ in a connection component containing an odd cycle are forced to be zero by equations \eqref{Sigma_sys}. The elements in a bipartite connection component must have equal absolute values, but the elements in the two partition classes of the component have opposite signs. Hence the value of an arbitrary element in the bipartite connection component can be chosen at will, while all other elements in the component are determined by the value of that first one.

By virtue of Theorem \ref{th_extremal} extremality of $\Sigma$ entails that system \eqref{Sigma_sys} has a one-dimensional solution space, namely the multiples of the matrix $\Sigma$ itself. Hence the graph $G$ has exactly one bipartite connection component. A solution $X$ of system \eqref{Sigma_sys}, in particular the matrix $\Sigma$ itself, must then have zero entries at all vertices which are not in this bipartite component, and all remaining entries have equal absolute value. Since $\Sigma$ has a unit diagonal, its entries can therefore assume only the values $-1,0,1$. This completes the proof.
\end{proof}

We are now able to formulate the main result.

{\theorem Let $A \in {\cal C}^n$ be an extremal copositive matrix. Then the following are equivalent:
\begin{itemize}
\item[(i)] All minimal zeros of $A$ have support of cardinality two.
\item[(ii)] There exists a positive definite diagonal matrix $D$ and an extremal copositive matrix $\Sigma$ with unit diagonal and all elements from the set $\{-1,0,1\}$ such that $A = D\Sigma D$.
    \end{itemize} }
    
\begin{proof}
The implication $(i) \Rightarrow (ii)$ is the assertion of Lemma \ref{1to2}.

Assume condition $(ii)$. By Lemma \ref{2to1} all minimal zeros of $\Sigma$ have supports of cardinality two. But then the same holds for $A$, because the minimal zero support set is preserved by automorphisms of ${\cal C}^n$ of the form $X \mapsto DXD$. This proves $(i)$.
\end{proof}

This result completely characterizes the class of extremal copositive matrices with minimal zero supports of cardinality two. Namely, these matrices are diagonally scaled versions of the extremal $\{-1,0,1\}$-matrices with unit diagonal, which have been already classified in~\cite{HoffmanPereira73}.

\bibliographystyle{plain}
\bibliography{copositive}

\end{document}